\documentclass{article}%
\usepackage{amsfonts}
\usepackage{amssymb}
\usepackage{graphicx}
\usepackage{amsmath}%
\newtheorem{theorem}{Theorem}

\newtheorem{definition}[theorem]{Definition}

\begin{document}

\title{Interlace polynomials and Tutte polynomials}
\author{Lorenzo Traldi\\Lafayette College\\Easton, Pennsylvania 18042}
\date{}
\maketitle

\begin{abstract}
Let $G$ be a graph with adjacency matrix $A(G)$. Consider the matrix
$IA(G)=\left(  I\text{ }|\text{ }A(G)\right)  $, where $I$ is the identity
matrix, and let $M(IA(G))$ be the binary matroid represented by $IA(G)$. Then
suitable parametrized versions of the Tutte polynomial of $M(IA(G))$ yield the
interlace polynomials of $G$, introduced by Arratia, Bollob\'{a}s and Sorkin
[J. Combin. Theory Ser. B 92 (2004) 199-233; Combinatorica 24 (2004) 567-584].
Interlace polynomials subsequently introduced by other authors may be obtained
from parametrized Tutte polynomials of the binary matroid represented by
$\left(  I\text{ }|\text{ }A(G)\text{ }|\text{ }I+A(G)\right)  $.

\bigskip


\bigskip

Mathematics Subject\ Classification. 05C50

\end{abstract}

\section{Introduction}

Motivated by problems that arise in the study of DNA sequencing, Arratia,
Bollob\'{a}s and Sorkin introduced a one-variable graph polynomial, the
\emph{vertex-nullity interlace polynomial}, in \cite{A1}. In subsequent work
\cite{A2, A} they observed that this one-variable polynomial may be obtained
from the Tutte-Martin polynomial of isotropic systems studied by\ Bouchet
\cite{Bi3}, introduced an extended two-variable version of the interlace
polynomial, and observed that the interlace polynomials are given by formulas
that involve the nullities of matrices over the two-element field, $GF(2)$.
Inspired by these ideas, Aigner and van der Holst \cite{AH}, Courcelle
\cite{C} and the author \cite{Tw, Tint} introduced several different
variations on the interlace polynomial theme.

All these references share the underlying presumption that although the theory
of the interlace polynomials is connected to that of the Tutte polynomial in
some ways, the two theories are largely separate in general. In this short
note we point out that in fact, the interlace polynomials of graphs can be
derived from parametrized Tutte polynomials of binary matroids associated with
adjacency matrices. We presume the reader is familiar with the standard
terminology of graph theory and matroid theory; see \cite{GM, O, We, W1} for instance.

\section{The identity-adjacency matroid}

We restrict our attention to looped simple graphs. That is, a graph $G$ is
given by specifying a finite set $V(G)$ of \emph{vertices}, declaring that
certain vertices are \emph{looped} and the others are not, and declaring that
certain pairs of distinct vertices are \emph{neighbors} and the other pairs
are not. The \emph{adjacency matrix }of $G$ is the $V(G)\times V(G)$ matrix
$A(G)$ with entries in $GF(2)$ given by: a diagonal entry is 1 if and only if
the corresponding vertex is looped, and an off-diagonal entry is 1 if and only
if the corresponding vertices are neighbors.

If $I$ is the $\left\vert V(G)\right\vert \times\left\vert V(G)\right\vert $
identity matrix then $IA(G)$ is the $\left\vert V(G)\right\vert \times
2\left\vert V(G)\right\vert $ matrix
\[
IA(G)=\left(  I\text{ }|\text{ }A(G)\right)  \text{.}%
\]
For convenience of notation and in order to indicate the relationship with our
previous work on the interlace polynomials, we use the Greek letters $\phi$
and $\chi$ to refer to the columns of the indicated submatrices of $IA(G)$:
the column of $I$ corresponding to $v$ is denoted $v_{\phi}$, and the column
of $A(G)$ corresponding to $v$ is denoted $v_{\chi}$.

\begin{definition}
The \emph{identity-adjacency matroid} $M(IA(G))$ is the binary matroid
represented by $IA(G)$.
\end{definition}

That is, $M(IA(G))$ is a matroid on the ground set $W(G)=\{v_{\phi}$,
$v_{\chi}\mid v\in V(G)\}$, and if $T\subseteq W(G)$ then the rank $r^{G}(T)$
of $T$ in $M(IA(G))$ equals the dimension of the $GF(2)$-vector space spanned
by the columns of $IA(G)$ corresponding to elements of $T$.

One way to define the \emph{Tutte polynomial} of $M(IA(G))$ is a polynomial in
the variables $s$ and $z$, given by the subset expansion%
\[
t(M(IA(G))=\sum_{T\subseteq W(G)}s^{r^{G}(W(G))-r^{G}(T)}z^{\left\vert
T\right\vert -r^{G}(T)}\text{.}%
\]
We do not give a general account of this famous invariant of graphs and
matroids here; thorough introductions may be found in \cite{Bo, BO, D, GM}.

Tutte polynomials of graphs and matroids are remarkable both for the amount of
structural information they contain and for the range of applications in which
they appear. Some applications (electrical circuits, knot theory, network
reliability, and statistical mechanics, for instance) involve graphs or
networks whose vertices or edges have special attributes of some kind --
impedances and resistances in circuits, crossing types in knot diagrams,
probabilities of failure and successful operation in reliability, bond
strengths in statistical mechanics. A natural way to think of these attributes
is to allow each element to carry two parameters, $a$ and $b$ say, with $a$
contributing to the terms of the Tutte polynomial corresponding to subsets
that include the given element, and $b$ contributing to the terms of the Tutte
polynomial corresponding to subsets that do not. Zaslavsky \cite{Z} calls the
resulting polynomial%
\begin{equation}
\sum_{T\subseteq W(G)}\left(
{\displaystyle\prod\limits_{t\in T}}
a(t)\right)  \left(
{\displaystyle\prod\limits_{w\notin T}}
b(w)\right)  s^{r^{G}(W(G))-r^{G}(T)}z^{\left\vert T\right\vert -r^{G}(T)}
\label{partutte}%
\end{equation}
the \emph{parametrized rank polynomial} of $M(IA(G))$; we denote it
$\tau(M(IA(G))$.

We do not give a general account of the theory of parametrized Tutte
polynomials here; the interested reader is referred to the literature, for
instance \cite{BR, EMT, Sok, T, Z}. However it is worth taking a moment to
observe that parametrized polynomials are very flexible, and the same
information can be formulated in many ways.\ For instance if $s$ and the
parameter values $b(w)$ are all invertible then formula (\ref{partutte}) is
equivalent to
\[
s^{r^{G}(W(G))}\cdot\left(
{\displaystyle\prod\limits_{w\in W(G)}}
b(w)\right)  \cdot\sum_{T\subseteq W(G)}\left(
{\displaystyle\prod\limits_{t\in T}}
\left(  \frac{a(t)}{b(t)s}\right)  \right)  (sz)^{\left\vert T\right\vert
-r^{G}(T)}\text{,}%
\]
which expresses $\tau(M(IA(G))$ as the product of a prefactor and a sum that
is essentially a parametrized rank polynomial with only $a$ parameters and one
variable, $sz$. We prefer formula (\ref{partutte}), though, because we do not
want to assume invertibility of the $b$ parameters.

Suppose that the various parameter values $a(w)$ and $b(w)$ are independent
indeterminates, and let $P$ denote the ring of polynomials with integer
coefficients in the set of $2+4\left\vert V(G)\right\vert $ independent
indeterminates $\{s,z\}\cup\{a(w)$, $b(w)\mid w\in W(G)\}$. Let $J$ be the
ideal of $P$ generated by the set of $2\left\vert V(G)\right\vert $ products
$\{a(v_{\phi})a(v_{\chi})$, $b(v_{\phi})b(v_{\chi})\mid v\in V(G)\}$, and let
$\pi:P\rightarrow P/J$ be the canonical map onto the quotient. Then the only
summands of (\ref{partutte}) that make nonzero contributions to $\pi
\tau(M(IA(G))$ correspond to subsets $T\subseteq W(G)$ with the property that
$\left\vert T\cap\{v_{\phi},v_{\chi}\}\right\vert =1$ $\forall v\in V(G)$.
Each such $T$ is a transversal of the partition of $W(G)$ into 2-element
subsets $\{v_{\phi},v_{\chi}\}$ corresponding to vertices of $G$; we denote
the collection of all such transversals $\mathcal{T}(W(G))$. Each
$T\in\mathcal{T}(W(G))$ has $\left\vert T\right\vert =\left\vert
V(G)\right\vert =r^{G}(W(G))$, so $s$ and $z$ have the same exponent in the
corresponding term of $\pi\tau(M(IA(G))$:
\[
\pi\tau(M(IA(G))=\pi\left(  \sum_{T\in\mathcal{T}(W(G))}\left(
{\displaystyle\prod\limits_{t\in T}}
a(t)\right)  \left(
{\displaystyle\prod\limits_{w\notin T}}
b(w)\right)  (sz)^{\left\vert V(G)\right\vert -r^{G}(T)}\right)  \text{.}%
\]

Every generator of $J$ involves the product of two different parameters
$a(v_{i})$, $b(v_{i})$ corresponding to a single vertex $v$ of $G$. It follows
that $\pi$ is injective when restricted to the additive subgroup $A$ of $P$
generated by products%
\[
\left(
{\displaystyle\prod\limits_{t\in T}}
a(t)\right)  \left(
{\displaystyle\prod\limits_{w\notin T}}
b(w)\right)  (sz)^{k}%
\]
where $k\geq0$ and $T\in\mathcal{T}(W(G))$. Consequently there is a
well-defined isomorphism of abelian groups $\pi^{-1}:\pi(A)\rightarrow A$, and
we have
\begin{equation}
\pi^{-1}\pi\tau(M(IA(G))=\sum_{T\in\mathcal{T}(W(G))}\left(
{\displaystyle\prod\limits_{t\in T}}
a(t)\right)  \left(
{\displaystyle\prod\limits_{w\notin T}}
b(w)\right)  (sz)^{\left\vert V(G)\right\vert -r^{G}(T)}\text{.}
\label{partutte2}%
\end{equation}

Note that $\pi^{-1}\pi\tau(M(IA(G))$, the image of the parametrized Tutte
polynomial $\tau(M(IA(G))$ under the mappings $\pi$\ and $\pi^{-1}$, might
also be described as the \emph{section} of $\tau(M(IA(G))$ corresponding to
$\mathcal{T}(W(G))$. Either way, formula (\ref{partutte2}) describes an
element of $P$, where $s$, $z$ and the various parameter values $a(w)$, $b(w)$
are all independent indeterminates.

Arratia, Bollob\'{a}s and Sorkin \cite{A} define the two-variable
\emph{interlace polynomial }$q(G)$ by the formula%
\begin{align*}
q(G)  &  =\sum_{S\subseteq V(G)}\left(  x-1\right)  ^{r(A(G)[S])}\left(
y-1\right)  ^{\left\vert S\right\vert -r(A(G)[S])}\\
&  =\sum_{S\subseteq V(G)}\left(  \frac{y-1}{x-1}\right)  ^{\left\vert
S\right\vert -r(A(G)[S])}\left(  x-1\right)  ^{\left\vert S\right\vert
}\text{.}%
\end{align*}
Here $r(A(G)[S])$ denotes the $GF(2)$-rank of the principal submatrix of
$A(G)$ involving rows and columns corresponding to vertices from $S$.

For $T\in\mathcal{T}(W(G))$ let $S(T)=\{v\in V(G)\mid v_{\chi}\in T\}$; then
$T\mapsto S(T)$ defines a bijection from $\mathcal{T}(W(G))$ onto the
power-set of $V(G)$. As $r^{G}(T)$ is the $GF(2)$-rank of the matrix%
\[
\left(  \text{columns }v_{\phi}\text{ with }v\notin S(T)\mid\text{columns
}v_{\chi}\text{ with }v\in S(T)\right)
\]
and the columns $v_{\phi}$ are columns of the identity matrix,
\[
r^{G}(T)=\left\vert V(G)\right\vert -\left\vert S(T)\right\vert
+r(A(G)[S(T)]).
\]
It follows that $q(G)$ may be obtained from $\pi^{-1}\pi\tau(M(IA(G))$ by
setting $a(v_{\phi})\equiv1$, $a(v_{\chi})\equiv x-1$, $b(v_{\phi})\equiv1$,
$b(v_{\chi})\equiv1$, $s=y-1$ and $z=1/(x-1)$. These assignments are not
unique; for instance the values of $s$ and $z$ may be replaced by
$s=(y-1)/\sigma$ and $z=\sigma/(x-1)$ for any invertible $\sigma$.

\section{Recursive formulas}

In this section we show how a recursive description of the parametrized rank
polynomial yields the recursive description of the interlace polynomial given
by Arratia, Bollob\'{a}s and Sorkin \cite{A}. Recall that a \emph{coloop} of a
matroid is an element that is included in every basis, and a \emph{loop} is an
element that is excluded from every basis. Suppose $r^{M}$ is the rank
function of a matroid $M$ on the ground set $W$. If $w\in W$ then $M-w$ is the
matroid on $W-\{w\}$ whose rank function is given by $r^{M}(T)$, for
$T\subseteq W-\{w\}$; $M/w$ is the matroid on $W-\{w\}$ whose rank function is
given by $r^{M}(T\cup\{w\})-r^{M}(\{w\})$.

Parametrized rank polynomials may be calculated recursively as follows:

\begin{enumerate}
\item If $\varnothing$ is the empty matroid then $\tau(\varnothing)=1.$

\item If $M$ is a matroid on $W$ and $w$ is a coloop of $M$ then
$\tau(M)=(a(w)+sb(w))\cdot\tau(M/w)$.

\item If $w$ is a loop of $M$ then $\tau(M)=(b(w)+za(w))\cdot\tau(M-w)$.

\item If $w$ is neither a coloop nor a loop of $M$ then $\tau(M)=b(w)\tau
(M-w)+a(w)\tau(M/w)$.
\end{enumerate}

Suppose $G$ is a\ graph and $v\in V(G)$ is not isolated; let $w$ be a neighbor
of $v$. If $v$ is looped then $\{x_{\phi}\mid x\in V(G)\}$ and $\{v_{\chi
}\}\cup\{x_{\phi}\mid x\neq v\}$ are both bases of $M(IA(G))$, because $I$
and
\[%
\begin{pmatrix}
1 & \mathbf{0}\\
\mathbf{\ast} & I^{\prime}%
\end{pmatrix}
\]
are both of rank $\left\vert V(G)\right\vert $, where $I^{\prime}$ is the
identity matrix of order $\left\vert V(G)\right\vert -1$. If $v$ is unlooped
then $\{x_{\phi}\mid x\in V(G)\}$ and $\{v_{\chi},w_{\chi}\}\cup\{x_{\phi}\mid
v\neq x\neq w\}$ are both bases of $M(IA(G))$, because $I$ and
\[%
\begin{pmatrix}
0 & 1 & 0\\
1 & \ast & 0\\
\mathbf{\ast} & \ast & I^{\prime\prime}%
\end{pmatrix}
\text{ }%
\]
are both of rank $\left\vert V(G)\right\vert $, where $I^{\prime\prime}$ is
the identity matrix of order $\left\vert V(G)\right\vert -2$. In either case,
we see that $v_{\phi}$ is not a coloop or a loop of $M(IA(G))$.

If $v$ is looped then $\{v_{\chi}\}\cup\{x_{\phi}\mid x\neq v\}$ and
$\{w_{\chi}\}\cup\{x_{\phi}\mid x\neq v\}$ are both bases of $M(IA(G))-v_{\phi
}$. If $v$ is not looped then $\{w_{\chi}\}\cup\{x_{\phi}\mid x\neq v\}$ and
$\{v_{\chi},w_{\chi}\}\cup\{x_{\phi}\mid v\neq x\neq w\}$ are both bases of
$M(IA(G))-v_{\phi}$. In either case, we see that $v_{\chi}$ is not a coloop or
a loop of $M(IA(G))-v_{\phi}$.

As $v_{\phi}$ has only one nonzero entry, the definition of matroid
contraction mentioned at the beginning of this section tells us that the rank
function of $M(IA(G))/v_{\phi}$ is the function on $W(G)-\{v_{\phi}\}$ defined
using the columns of the matrix $IA(G)^{\prime}$ obtained from $IA(G)$ by
removing both the column $v_{\phi}$ and the row corresponding to $v$. In
particular, the rank of the whole matroid is $\left\vert V(G)\right\vert -1$.
As $\{x_{\phi}\mid x\neq v\}$ and $\{v_{\chi}\}\cup\{x_{\phi}\mid v\neq x\neq
w\}$ are both bases, $v_{\chi}$ is not a coloop or a loop of $M(IA(G))/v_{\phi
}$.

Let $G-v$ be the graph obtained from $G$ by removing $v$ and all edges
incident on it. Then $IA(G-v)$ is the matrix obtained by removing the column
$v_{\chi}$ from $IA(G)^{\prime}$, so
\[
M(IA(G-v))=(M(IA(G))/v_{\phi})-v_{\chi}.
\]
Using step 4 of the recursion to remove $v_{\phi}$ and then $v_{\chi}$, we see
that%
\begin{gather*}
\pi\tau(M(IA(G)))=\\
\pi a(v_{\phi})b(v_{\chi})\tau((M(IA(G))/v_{\phi})-v_{\chi})+\pi b(v_{\phi
})a(v_{\chi})\tau((M(IA(G))-v_{\phi})/v_{\chi})\\
=\pi a(v_{\phi})b(v_{\chi})\tau(M(IA(G-v)))+\pi b(v_{\phi})a(v_{\chi}%
)\tau((M(IA(G))-v_{\phi})/v_{\chi})\text{.}%
\end{gather*}
With the parameter values given at the end of Section 2, this yields%
\[
\pi^{-1}\pi\tau(M(IA(G)))=\pi^{-1}\pi\tau(M(IA(G-v)))+(x-1)\pi^{-1}\pi
\tau((M(IA(G))-v_{\phi})/v_{\chi})
\]
or equivalently,%
\begin{equation}
q(G)=q(G-v)+(x-1)\pi^{-1}\pi\tau((M(IA(G))-v_{\phi})/v_{\chi}).
\label{recurs1}%
\end{equation}

Suppose $v$ is looped, and let $G^{v}$ be the graph obtained from $G$ by
toggling all adjacencies between neighbors of $v$, and also toggling the loop
status of every neighbor of $v$. Consider two matrices%
\[%
\begin{pmatrix}
1 & \mathbf{1} & \mathbf{0}\\
\mathbf{1} & B & C\\
\mathbf{0} & D & E
\end{pmatrix}
\text{ and }%
\begin{pmatrix}
1 & \mathbf{0} & \mathbf{0}\\
\mathbf{1} & \overline{B} & C\\
\mathbf{0} & D & E
\end{pmatrix}
\]
where bold numerals indicate row and column vectors and the overbar indicates
a submatrix in which all entries have been toggled (reversed). Elementary
column operations tell us that the two matrices have the same $GF(2)$-rank.
Consequently, if $v_{\chi}\in T\in\mathcal{T}(W(G))$ then $r^{G}%
(T)=1+r^{G^{v}-v}(T-\{v_{\chi}\})$. The definition of matroid contraction
tells us that if $v_{\chi}\in T\in\mathcal{T}(W(G))$ then the rank of
$T-\{v_{\chi}\}$ in $(M(IA(G))-v_{\phi})/v_{\chi}$ is $r^{G}(T)-1$. It follows
that the ranks of $T-\{v_{\chi}\}$ in $(M(IA(G))-v_{\phi})/v_{\chi}$ and
$M(IA(G^{v}-v))$ are equal. Combining this equality with (\ref{recurs1}), we
see that if $v$ is looped then
\begin{align}
q(G)  &  =q(G-v)+(x-1)\pi^{-1}\pi\tau(M(IA(G^{v}-v)))\label{recurs2}\\
&  =q(G-v)+(x-1)q(G^{v}-v).\nonumber
\end{align}

Formula (\ref{recurs2}) is one of the two fundamental recursive formulas for
$q(G)$ \cite{A}.

Deriving the other fundamental recursive formula takes a little more work,
because the term $\pi^{-1}\pi\tau((M(IA(G))-v_{\phi})/v_{\chi})$ in formula
(\ref{recurs1}) does not correspond to a single interlace polynomial.

Suppose that $v$ is unlooped and has an unlooped neighbor $w$. As in \cite{A},
say that two vertices $x,y\notin\{v,w\}$ are \emph{distinguished} by $\{v,w\}$
if they have different, nonempty neighborhoods in $\{v,w\}$, and let $G^{vw}$
denote the graph obtained from $G$ by toggling all adjacencies between
vertices distinguished by $\{v,w\}$.

The matroid $(M(IA(G))-v_{\phi})/v_{\chi}$ has $W(G)-\{v_{\phi},v_{\chi
}\}=W(G-v)$ as its ground set, and $\pi^{-1}\pi\tau((M(IA(G))-v_{\phi
})/v_{\chi})$ includes nonzero contributions from elements of $\mathcal{T}%
(G-v)$. Split $\pi^{-1}\pi\tau((M(IA(G))-v_{\phi})/v_{\chi})$ into two parts,
$S_{\phi}$ and $S_{\chi}$, with $S_{\phi}$ including the contributions from
elements of $\mathcal{T}(G-v)$ that include $w_{\phi}$ and $S_{\chi}$
including the contributions from elements of $\mathcal{T}(G-v)$ that include
$w_{\chi}$.

Consider three matrices%
\begin{gather*}
\left(
\begin{array}
[c]{cccccc}%
0 & 0 & \mathbf{1} & \mathbf{1} & \mathbf{0} & \mathbf{0}\\
1 & 1 & \mathbf{1} & \mathbf{0} & \mathbf{1} & \mathbf{0}\\
\mathbf{1} & \mathbf{0} & B_{11} & B_{12} & B_{13} & B_{14}\\
\mathbf{1} & \mathbf{0} & B_{21} & B_{22} & B_{23} & B_{24}\\
\mathbf{0} & \mathbf{0} & B_{31} & B_{32} & B_{33} & B_{34}\\
\mathbf{0} & \mathbf{0} & B_{41} & B_{42} & B_{43} & B_{44}%
\end{array}
\right)  \text{, }\left(
\begin{array}
[c]{cccccc}%
0 & 0 & \mathbf{1} & \mathbf{1} & \mathbf{0} & \mathbf{0}\\
1 & 1 & \mathbf{0} & \mathbf{0} & \mathbf{0} & \mathbf{0}\\
\mathbf{1} & \mathbf{0} & \overline{B}_{11} & B_{12} & \overline{B}_{13} &
B_{14}\\
\mathbf{1} & \mathbf{0} & \overline{B}_{21} & B_{22} & \overline{B}_{23} &
B_{24}\\
\mathbf{0} & \mathbf{0} & B_{31} & B_{32} & B_{33} & B_{34}\\
\mathbf{0} & \mathbf{0} & B_{41} & B_{42} & B_{43} & B_{44}%
\end{array}
\right) \\
\text{and }\left(
\begin{array}
[c]{cccccc}%
0 & 0 & \mathbf{1} & \mathbf{1} & \mathbf{0} & \mathbf{0}\\
1 & 1 & \mathbf{0} & \mathbf{0} & \mathbf{0} & \mathbf{0}\\
\mathbf{1} & \mathbf{0} & B_{11} & \overline{B}_{12} & \overline{B}_{13} &
B_{14}\\
\mathbf{1} & \mathbf{0} & \overline{B}_{21} & B_{22} & \overline{B}_{23} &
B_{24}\\
\mathbf{0} & \mathbf{0} & \overline{B}_{31} & \overline{B}_{32} & B_{33} &
B_{34}\\
\mathbf{0} & \mathbf{0} & B_{41} & B_{42} & B_{43} & B_{44}%
\end{array}
\right)  \text{ .}%
\end{gather*}
Elementary column operations tell us that the first two matrices have the same
$GF(2)$-rank, and elementary row operations tell us that the third also has
the same $GF(2)$-rank. Removing the second row and column from the third
matrix reduces the $GF(2)$-rank by 1, clearly. It follows that if $v_{\chi
},w_{\phi}\in T\in\mathcal{T}(W(G))$ then%
\[
r^{G}(T)-1=r^{G^{vw}-w}(T-\{w_{\phi}\})\text{.}%
\]

Recall that $r^{G}(T)-1$ is the rank of $T-\{v_{\chi}\}$ in $(M(IA(G))-v_{\phi
})/v_{\chi}$, and the parameter values given in Section 2 include $a(v_{\chi
})=x-1$ and $a(w_{\phi})=1$. It follows that the contribution of $T-\{w_{\phi
}\}$ to $\pi^{-1}\pi\tau(M(IA(G^{vw}-w)))$ is the product of $x-1$ and the
contribution of $T-\{v_{\chi}\}$ to $\pi^{-1}\pi\tau((M(IA(G))-v_{\phi
})/v_{\chi})$. Consequently, if we split $q(G^{vw}-w)$ into two parts,
$q_{\phi}$ and $q_{\chi}$, with $q_{\phi}$ including the contributions from
elements of $\mathcal{T}(G^{vw}-w)$ that include $v_{\phi}$ and $q_{\chi}$
including the contributions from elements of $\mathcal{T}(G^{vw}-w)$ that
include $v_{\chi}$, then
\[
(x-1)S_{\phi}=q_{\chi}.
\]
If $v_{\phi}\in T\in\mathcal{T}(G^{vw}-w)$, then the corresponding column has
only one nonzero entry; elementary column operations show that the rank of $T$
in $M(IA(G^{vw}-w))$ is 1 more than the rank of $T-\{v_{\phi}\}$ in
$M(IA(G^{vw}-v-w))$. It follows that
\[
q_{\phi}=q(G^{vw}-v-w)\text{.}%
\]
Formula (\ref{recurs1}) now tells us that
\begin{align}
q(G)  &  =q(G-v)+(x-1)(S_{\phi}+S_{\chi})=q(G-v)+q_{\chi}+(x-1)S_{\chi
}\label{recurs3}\\
&  =q(G-v)+q(G^{vw}-w)-q_{\phi}+(x-1)S_{\chi}\nonumber\\
&  =q(G-v)+q(G^{vw}-w)-q(G^{vw}-v-w)+(x-1)S_{\chi}\text{.}\nonumber
\end{align}

It remains only to discuss $S_{\chi}$. Consider two matrices%
\[
\left(
\begin{array}
[c]{cccccc}%
0 & 1 & \mathbf{1} & \mathbf{1} & \mathbf{0} & \mathbf{0}\\
1 & 0 & \mathbf{1} & \mathbf{0} & \mathbf{1} & \mathbf{0}\\
\mathbf{1} & \mathbf{1} & B_{11} & B_{12} & B_{13} & B_{14}\\
\mathbf{1} & \mathbf{0} & B_{21} & B_{22} & B_{23} & B_{24}\\
\mathbf{0} & \mathbf{1} & B_{31} & B_{32} & B_{33} & B_{34}\\
\mathbf{0} & \mathbf{0} & B_{41} & B_{42} & B_{43} & B_{44}%
\end{array}
\right)  \text{ and }\left(
\begin{array}
[c]{cccccc}%
0 & 1 & \mathbf{0} & \mathbf{0} & \mathbf{0} & \mathbf{0}\\
1 & 0 & \mathbf{0} & \mathbf{0} & \mathbf{0} & \mathbf{0}\\
\mathbf{1} & \mathbf{1} & B_{11} & \overline{B}_{12} & \overline{B}_{13} &
B_{14}\\
\mathbf{1} & \mathbf{0} & \overline{B}_{21} & B_{22} & \overline{B}_{23} &
B_{24}\\
\mathbf{0} & \mathbf{1} & \overline{B}_{31} & \overline{B}_{32} & B_{33} &
B_{34}\\
\mathbf{0} & \mathbf{0} & B_{41} & B_{42} & B_{43} & B_{44}%
\end{array}
\right)  \text{ .}%
\]
Elementary column operations tell us that the two matrices have the same
$GF(2)$-rank. It follows that if $v_{\chi},w_{\chi}\in T\in\mathcal{T}(W(G))$
then%
\[
r^{G}(T)=2+r^{G^{vw}-v-w}(T-\{v_{\chi},w_{\chi}\})\text{.}%
\]
The rank of $T-\{v_{\chi}\}$ in $(M(IA(G))-v_{\phi})/v_{\chi}$ is
\[
r^{G}(T)-1=1+r^{G^{vw}-v-w}(T-\{v_{\chi},w_{\chi}\})
\]
so considering the parameter values given at the end of Section 2, we see that
the contribution of $T-\{v_{\chi}\}$ to $S_{\chi}$ is the product of $x-1$ and
the contribution of $T-\{v_{\chi},w_{\chi}\}$ to $\pi^{-1}\pi\tau
(M(IA(G^{vw}-v-w))$. It follows that
\[
S_{\chi}=(x-1)q(G^{vw}-v-w)\text{,}%
\]
so formula (\ref{recurs3}) tells us that
\[
q(G)=q(G-v)+q(G^{vw}-w)-q(G^{vw}-v-w)+(x-1)^{2}q(G^{vw}-v-w)\text{.}%
\]
This is the second fundamental recursive formula for $q(G)$ given in \cite{A}.

\section{The identity-adjacency-sum matroid}

We call the binary matroid $M(IAS(G))$ represented by the matrix%
\[
IAS(G)=\left(  I\mid A(G)\mid A(G)+I\right)
\]
the \emph{identity-adjacency-sum matroid} of the graph $G$. The column of
$A(G)+I$ corresponding to a vertex $v$ is denoted $v_{\psi}$, and the other
columns are denoted $v_{\phi}$, $v_{\chi}$ as in\ Section 2; then $W^{\prime
}(G)=\{v_{\phi}$, $v_{\chi}$, $v_{\psi}\mid v\in V(G)\}$ is the ground set of
$M(IAS(G))$. Let $P^{\prime}$ be the ring of polynomials with integer
coefficients in the set of $2+6\left\vert V(G)\right\vert $ independent
indeterminates $\{s,z\}\cup\{a(w)$, $b(w)\mid w\in W^{\prime}(G)\}$, let
$J^{\prime}$ be the ideal of $P^{\prime}$ generated by the set of $4\left\vert
V(G)\right\vert $ products $\{a(v_{\phi})a(v_{\chi})$, $a(v_{\phi})a(v_{\psi
})$, $a(v_{\chi})a(v_{\psi})$, $b(v_{\phi})b(v_{\chi})b(v_{\psi})\mid v\in
V(G)\}$, and let $\pi^{\prime}:P^{\prime}\rightarrow P^{\prime}/J^{\prime}$ be
the canonical map onto the quotient. Then the discussion of Section 2 is
readily modified to show that $(\pi^{\prime})^{-1}\pi^{\prime}\tau(M(IAS(G)))$
consists only of terms associated with
\[
\mathcal{T}^{\prime}(W^{\prime}(G))=\{T\subseteq W^{\prime}(G)\mid\left\vert
T\cap\{v_{\phi},v_{\chi},v_{\psi}\}\right\vert =1\text{ }\forall v\in
V(G)\text{.}%
\]

The reader familiar with the interlace polynomials introduced by Aigner and
van der Holst \cite{AH}, Courcelle \cite{C}, and the author \cite{Tint} will
have no trouble modifying the discussion of Section 2 to show that appropriate
values for $s$, $z$ and the $a$ and $b$ parameters yield all of these
interlace polynomials from the parametrized rank polynomial $\tau(M(IAS(G)))$.
Notice also that $M(IA(G))$ is a sub-matroid of $M(IAS(G)$, and $(\pi^{\prime
})^{-1}\pi^{\prime}\tau(M(IAS(G)))$ yields $\pi^{-1}\pi\tau(M(IA(G)))$ by
assigning $a(v_{\psi})\equiv0$ and $b(v_{\psi})\equiv1$; consequently the
theory given in Section 2 is contained in the one described here. (It is for
for expository convenience, not generality, that we detail the theory of
$M(IA(G))$ rather than that of $M(IAS(G))$.)

\section{Comments}

1. As explained by Ellis-Monaghan and Sarmiento \cite{EMS}, results of Las
Vergnas \cite{L2, L1} and Martin \cite{Ma} on circuit partitions of planar
4-regular graphs imply that if $G$ happens to be a circle graph obtained from
a planar 4-regular graph, then the vertex-nullity interlace polynomial of $G$
may be obtained from the \textquotedblleft diagonal\textquotedblright\ Tutte
polynomial of an associated checkerboard graph. (The diagonal\ Tutte
polynomial is obtained by setting the two variables of the ordinary
(non-parametrized) Tutte polynomial equal to each other.) This connection
cannot extend directly to non-planar graphs, as Martin \cite[p.76]{Ma} pointed
out, because the complete graph $K_{5}$ has too many Euler circuits to be
represented in a 5-element matroid.

Also, Aigner and van der Holst \cite{AH} observed that the vertex-nullity
interlace polynomial of a bipartite graph may be obtained from the diagonal
Tutte polynomial of an associated binary matroid. This result is connected to
the preceding paragraph through de Fraysseix's theorem connecting bipartite
circle graphs to planar graphs \cite{F}. More recently, Brijder and Hoogeboom
\cite{BH1} have introduced interlace polynomials for delta-matroids, and in
particular for matroids. They use them to extend the connection between
vertex-nullity interlace polynomials and diagonal Tutte polynomials to
arbitrary matroids.

One way to summarize the content of the present note is this: Using parameters
has the effect of algebraically restricting $\tau(M(IA(G))$ to $\mathcal{T}%
(W(G))$ and $\tau(M(IAS(G))$ to $\mathcal{T}^{\prime}(W^{\prime}(G))$, and
unlike restriction to the diagonal\ Tutte polynomial, these restrictions are
effective for all graphs and all interlace polynomials.

2. The matrix $IA(G)$ appears in \cite{AH} and \cite{A2}, together with the
observation that the vertex-nullity interlace polynomial of $G$ is equal to
Bouchet's \emph{Tutte-Martin polynomial} of the isotropic system associated
with the row space of $IA(G)$ \cite{Bi3, B5}.

The matrix $IAS(G)$ also appears in \cite{AH}, where Aigner and van der Holst
showed that their interlace polynomial $Q$ may be obtained by summing over
submatrices of $IAS(G)$ associated with elements of $\mathcal{T}^{\prime
}(W^{\prime}(G))$. The content of this note came to mind after R. Brijder
pointed out the appearance of the same submatrices in our work with
nonsymmetric modified interlacement matrices \cite{Tnew}, when we read
Bouchet's comment \cite{B1} that Eulerian multimatroids are \textquotedblleft
sheltered\textquotedblright\ by matroids and wondered whether $M(IAS(G))$ is
in general an appropriate \textquotedblleft sheltering\textquotedblright%
\ matroid for the 3-matroid associated with an isotropic system with
fundamental graph $G$. Indeed it is!

3. Some known properties of interlace polynomials can be readily explained
using known properties of Tutte polynomials. For instance, the analogy between
pendant-twin reductions for the interlace polynomial and series-parallel
reductions for the Tutte polynomial noted by Bl\"{a}ser and Hoffman \cite{BH},
Ellis-Monaghan and Sarmiento \cite{EMS} and the author \cite{Tw, Tint} is more
than an analogy: when $v$ is pendant on $w$, $v_{\chi}$ and $w_{\phi}$ are
parallel in $M(IA(G))$; and when $v$ and $w$ are twins, $v_{\chi}$ and
$w_{\chi}$ are parallel in $M(IA(G))/\{v_{\phi},w_{\phi}\}$. Known properties
of Tutte polynomials\ also provide new insights into interlace polynomials;
for instance, the interlace polynomials of $G$ have activities expansions with
respect to bases of $M(IA(G))$ or $M(IAS(G))$.

4. In closing: understanding more about the graph-theoretic significance of
the matroids $M(IA(G))$ and $M(IAS(G))$ would help in understanding the
significance of the interlace polynomials.

\end{document}